\documentclass[a4paper,twoside,12pt]{article}

\usepackage[english]{babel}

\usepackage{amsmath,amssymb}

\usepackage{amsthm}

\begin{document}
\title{Alternating groups and rational functions on surfaces.}
\author{ Sonia Brivio and Gian Pietro Pirola} 
\footnotetext[1]{Partially
supported by 1)Cofin 2003: {\em Spazi di moduli e teoria di Lie.} (Murst); 2) Gnsaga; 3) Far 2002 (Pavia): {\em
Variet\`{a} algebriche, calcolo algebrico, grafi orientati e topologici}.
}
\newtheorem{teo}{Theorem}[subsection]
\newtheorem{prop}[teo]{Proposition}
\newtheorem{lem}[teo]{Lemma}
\newtheorem{rem}[teo]{Remark}
\newtheorem{defi}[teo]{Definition}
\newtheorem{ex}[teo]{Example}

\maketitle
\begin{abstract}
Let $X$ be a smooth complex projective surface and let ${\bf C}(X)$ denote the field of rational functions on $X$. 
In this paper, we prove that for any $m > M(X)$, there exists a rational dominant map $f \colon X \to Y$, which is 
generically finite of degree $m$, into a complex rational ruled surface $Y$, whose  monodromy  is the alternating group 
$A_m$. This gives a finite algebraic extension  ${\bf C}(X) \colon  {\bf C}(x_1,x_2)$ of degree $m$, whose normal closure 
has Galois group $A_m$.
 \end{abstract}
\section{Introduction. } Let  $F$ be an extension field of $L$, we denote by  $G(F \colon L)$ the Galois group  of the extension $F \colon L$, which   consists of all automorphisms of the field $F$ which fix $L$ elementwise.
If  $F \colon L$ is finite and separable, its normal closure $N \colon L$ is a Galois extension, see \cite{Garling}. Set $M(F , L) = G(N \colon L)$. Let  $X$ be  an irreducible complex algebraic variety, we can associate to it  the field ${\bf C}(X)$ of rational functions on $X$. This gives  a one to one correspondence between birational classes of irreducible complex algebraic varieties and finitely generated extensions of ${\bf C}$. Let $X$ and $Y$ be irreducible complex algebraic varieties of the same dimension $n$. Let $f \colon X \to Y$ be a generically finite dominant morphism of degree $d$. The  field  ${\bf C}(X)$ is a finite algebraic extension of degree $d$ of the field ${\bf C}(Y)$, the group $M({\bf C}(X) , {\bf C}(Y))$  is called the {\it Galois group of the morphism $f$ }, see \cite{Harris}. There is an   isomorphism between the Galois group of  $f$ and  the monodromy group $M(f )$, associated to the topological covering induced by $f$, see 2.1. Fix an irreducible variety $X$ of dimension $n$, ${\bf C}(X)$ is a finite algebraic extension of ${\bf C}({\bf P}^n)= {\bf C}(x_1, x_2,...,x_n)$, see \cite{Zariski}. The study of  possible monodromy groups for $X$   is a classic, algebraic and geometric  problem. In general, $M(f)$ is a subgroup of the symmetric group $S_d$. It is interesting to see in which cases  $M(f )$ is a subgroup  of the alternating group $A_d $; if this happens  we   say that   $f$  has {\it even monodromy}. 
\hfill \par
If $n =1$: let $X$ be a compact Riemann surface of genus $g $. Any  non constant meromorphic function $f \in  {\bf C}(X)$,  of degree $d$,   gives a holomorphic map $f \colon X \to {\bf P}^1$, which is a ramified covering of degree $d$.  $f$ is indecomposable if and only if the group $M(f)$ is a primitive subgroup of $S_d$.  There are several results on even monodromy of such maps: first of all by Riemann's existence theorem, $\forall g \geq 0$ and $\forall d \geq 2g +3$,  there are Riemann surfaces of genus $g$ admitting maps with monodromy group $A_d$, see \cite{Fried}. Actually, for a generic Riemann surface  $X$ of genus $g \geq 4$, for any indecomposable map the monodromy group is either $A_d$ or the symmetric group $S_d$, see \cite{Guralnick} and \cite{GuralnickM}. Finally, a generic compact Riemann surface of genus $1$ admits meromorphic functions   with  monodromy group $A_d$, for $d \geq 4$, see \cite{FriedKlassen}. This result has been recently generalized to any compact Riemann surface $X$ of genus $g$ for $d \geq 12g +4$, see \cite{ArtebaniPirola}. This implies that every extension field $F \colon {\bf C}$, with trascendence degree $1$,  can be realized as  a finite algebraic extension  of degree $d$, $F \colon L$, with   $L \simeq {\bf C}(x)$ and  monodromy group $ M(F,L) = A_d$.    \hfill \par 
In higher  dimension there are many various results concerning the   monodromy of branched coverings $f \colon X \to Y$  of a variety $ Y$ ( multiple planes theory, braid groups, Chisini problem, fundamental groups of the complement of a divisor ..., see \cite{Chisini}, \cite{Severi}, \cite{Catanese2}, \cite{Abhyankar}, \cite{Nori}, etc.). On the other hand, not much seems to be known when $X$ is fixed, for instance, the existence of maps $f$  with $M(f)$ solvable  is unknown also for projective surfaces of degree  $d \geq 6$. 
It is easy to produce, by general linear projections, finite maps $X \to {\bf P}^n$ with monodromy  the full symmetric group $S_d$, see \cite{Severi2}. So it is interesting to see if other primitive  groups can be realized as monodromy of $X$. In this paper, we  deal with  surfaces and even monodromy groups. Our  result is the following:
\hfill\par\hfill\par
{\bf Theorem 1} {\it  Let $F$ be an extension field of ${\bf C}$, with trascendence degree $2$. 
Then there exists an integer $ M(F)$ with the following property: for any $m > M(F)$,   
$F$ admits a subfield $L \simeq {\bf C}(x_1,x_2)$ such  that $F \colon L$ is a finite algebraic
 extension of degree $m$ and the group $M( F, L)$ is the alternating group $A_m$.}
\hfill \par \hfill\par
We will deduce theorem 1 from the following geometric result:
\hfill \par\hfill\par
{\bf Theorem 2} {\it Let $S$ be a smooth, complex, projective surface and $K_S$ denote a canonical divisor on $S$. 
Let $H$ be a very ample divisor on $S$, with $H^ 2 \geq 5$ and such that $(S, O_S(H))$ 
does not contain lines or conics. Set $g = p_a(2H + K_S)$. Then, for any $m > 16g +7$,  there exist a smooth complex 
projective surface $X$, in the birational class of $S$, and a generically finite surjective morphism 
$f \colon X \to Y$, of degree $m$, into a smooth complex rational ruled surface $Y$, such that the monodromy group
 $M(f)$ is the alternating group $A_m$. }
\hfill \par
Let us  describe briefly the method we  use in proving this result. Let $H$ be a very ample divisor on $S$: 
under our assumptions, which are actually verified by almost all $H$, we can find a Lefschetz pencil $P$ in the linear 
system  $ \vert 2H + K_S \vert $, whose elements are all irreducible, see 4.1.  By blowing up the base points of $P$, 
we produce a smooth, complex, projective surface $\hat S$,  and a surjective morphism $\phi \colon \hat S \to {\bf P}^1$,
 with fibre $F_t$, see 4.2. The pull back of $O_S(H + K_S)$ defines a natural spin bundle $L_t$ on each smooth fibre of 
$\phi$.  Following the method of \cite{ArtebaniPirola}, for each smooth fibre $F_t$, we can introduce the variety 
${\cal H}(F_t, D_t)$, parametrizing a family of meromorphic functions $f_t$ on $F_t$ with even monodromy, related to 
$L_t$, see 4.3. As $t$ varies on ${\bf P}^1$, we have a family $p \colon {\cal H} \to {\bf P}^1$ of projective varieties. 
  Our aim   is to glue these meromorphic functions in a suitable way. This can be done by  producing a section of $p$. 
Since for any $t \in {\bf P}^1$ the fibre $p^{-1}(t)$ is a  normal rationally connected variety, we can apply the 
following result: {\it every family of rationally  connected varieties over a smooth curve admits a section}. This property,
 conjectured by Koll$\acute a$r, has been recently proved by Graber, Harris, Starr, and by de Jong and Starr, 
(see \cite{GraberHarrisStarr} and  \cite{deJongStarr}). The existence of a section allows us to produce a generically finite 
surjective morphism $f \colon X \to Y$, where $X $ is birationally equivalent to  $S$, $Y$ is a smooth complex rational 
ruled surface,  such that the restriction $f_{\vert F_t}$ to a general  smooth fibre  has monodromy $A_m$. To conclude our
 proof, we show  that the monodromy of a general smooth fibre completely induces the monodromy of $f$. For this, we use  
a topological result of Nori, see 5.1. \hfill \par
Finally, we  apply our result to surfaces of general type  with ample canonical divisor,  see  5.3.  
We conjecture that theorem 1  holds for any finitely generated  extension $F$ of the complex field. 
\hfill \par
We would like to thank Enrico Schlesinger, who read the preliminary form of our manuscript and suggested many improvements, finally we are very grateful to the referee who pointed out that the proof of pr. 3.1.1 iii) was completely missing in the previous version of the paper.
\hfill\par 
\section{Preliminaries}
\subsection{Monodromy}
Let $X$ and $Y$ be irreducible complex algebraic varieties of the same dimension $n$. 
Let $f \colon X \to Y$ be a generically finite dominant morphism of degree $d$. We recall the definition of 
{\it monodromy group  $M (f)$ of $f$ }, see \cite{Harris}. Let $G $ be the  Galois group of the morphism $f$, see sec.1,
   $G$ acts faithfully on the general fibre ${f}^{-1}(y)$ and so can be seen as a subgroup of
 $Aut ( {f}^{-1}(y)) \simeq S_d$.   Let    $U \subset Y$ be an open dense subset such that the restriction
 $f \colon {f}^{-1}(U) \to U$ is a covering of degree $d$ in the classical topology (i.e. non ramified).
 For any point $y \in U$, let ${f}^{-1}(y) = \{ x_1,.., x_d \}$, we have the  monodromy representation of 
the fundamental group ${\pi}_1(U,y)$:
$$\rho(f,y) \colon {\pi}_1(U,y) \to Aut( {f}^{-1}(y)),$$
sending $[\alpha] \to \sigma(\alpha)$, where $\sigma(\alpha)$ is the  automorphism which sends $x_i$ to the end-point 
of the lift of $\alpha$ at the point $x_i$. Let $M(f,y) = \rho(f,y) ( {\pi}_1(U,y))$.  It is easy to verify that 
 $M(f,y)$ is isomorphic to the Galois group $G$ of $f$, and so does not depend on the choice of the open subset $U$.
The monodromy group $M(f )$  is defined as the coniugacy class of the transitive subgroups  $M(f,y)$. 
\subsection{Rational connectedness } Let $X$ be a  proper complex algebraic  variety of dimension $n$.  We recall that $X$ is  {\it rationally  connected} if and only if  for very general closed points $p$, $q \in X$ there is an 
 irreducible rational curve $C \subset X$ which contains $p$ and $q$, see \cite{Kollar}.
 \hfill\par
In the sequel we will need the following:
\hfill\par
\begin{prop} Let $X \subset{\bf P}^N$ be a complex irreducible variety of codimension $m$ which is the complete intersection 
of $ m$ hypersurfaces $Q_1$, $Q_2$,.., $Q_m$, of degree $d_1$, $d_2$,.., $d_m$. Let $ h = dim (Sing(X)) $  and $h = -1$ if $X$ is smooth. If 
\begin{equation}
   \sum_{i = 1}^{m}{d_i} \ + h + 1 \leq N 
\end{equation}
 then $X$ is rationally connected. In particular, a complete intersection $X$ of $m$ quadrics with $h \leq N -2m-1$ is rationally connected. 
\end{prop}
\begin{proof} 
If $X $ is   smooth, then  $\sum_{i = 1}^{m}{d_i} \leq N $ implies 
that  $X$ is a Fano variety, hence it is  rationally connected, (see \cite{Kollar} p.240.)
\hfill\par 
So we can assume that $dim(Sing(X)) = h \geq 0$. 
Let $H \in {({\bf P}^N)}^*$ be a  general hyperplane:  the hyperplane section $Y = X \cap H \subset {\bf P}^{N-1}$ is  a  complete  intersection, irreducible and   non degenerate,   of $m$ hypersurfaces of  ${\bf P}^{N-1}$ of degree $d_1$,..., $d_m$. Moreover, $dim (Sing (Y)) = h-1$ and $Y$ satisfies inequality $(1)$. Then it follows, by induction on $h$,  that $Y$  is   rationally connected. Finally,  for  general points $p$ and $q \in X$, there exists  a rationally  connected hyperplane section $Y$ containing  $p$ and $q$, hence an irreducible rational curve $C \subset X$ connecting the two points.
This concludes the proof.
\end{proof}
\hfill\par
\hfill\par
We remark that, in the provious proof, one can  intersect $X$ with a general linear space of  
dimension $N -h -1$ to get a smooth Fano variety connecting two general points of $X$.
 \hfill\par
An important  property of rational connectedness is given by the  following result,  see \cite{GraberHarrisStarr},  
\cite {deJongStarr}, and  \cite{KMM}:
\begin{teo} Let $p \colon X \to B$ be a proper flat morphism from a   complex projective variety into a smooth complex  
projective curve, assume that  $p$ is smooth over an open dense subset $U$ of $B$.  If the general fibre of $p$ is a normal
 and rationally connected variety, then $p$ has a section. Moreover, for any arbitrary finite set $A \subset U$ and for any
 section  ${\sigma}_1 \colon A \to p^{-1}(A)$,  there exists a section  $\sigma  \colon B \to X$ such that 
$\sigma_{\vert A} = {\sigma}_1$.
\end{teo}
\hfill\par
\subsection{Notations}  Let $S$ be a smooth, complex, connected, projective surface: we denote by $O_S$ the structure sheaf 
and  by $K_S$ a canonical divisor of $S$, so that $O_S(K_S) $ is the sheaf of the holomorphic $2$-forms. Let 
$q(S) = h^1(S,O_S)$ be the irregularity of $S$, let $p_g(S) = h^2(S,O_S)= h^0(S,O_S(K_S))$ be the geometric genus of $S$,
 finally let $p_n(S) = h^0(S,O_S(nK_S))$, $n \geq 1$, be the plurigenera of $S$. We denote by $k(S)$ the Kodaira dimension 
of $S$. A minimal surface $S$ is said {\it of general type } if  $k(S)= 2$. Let $C \subset S$ be an irreducible curve on $S$,
 we denote by $p_a(C) = h^1(C,O_C)$ the arithmetic genus of $C$, then $p_a(C) = 1 + \frac {1}{2}(C^2 + C \cdot K_S)$. If $C$ 
is smooth, $p_a(C) = g(C)$ is the geometric genus of $C$, and ${O_S(K_S+C)}_{\vert C} = {\omega}_C$ the canonical line 
bundle on $C$. A $g^r_d$ on a smooth curve $C$ is a linear serie (not necessarily  complete) on $C$ of degree $d$ and 
dimension exactly $r$. 
\subsection{Very ample line bundles} Let $L$ be a line bundle on $S$, $L$ is said $k$-spanned for $k \geq 0$ (i.e. it 
defines a $k$-th order embedding),  if for any distinct points $z_1,z_2,...z_t$ on $S$ and any positive integers 
$k_1,k_2,...k_t$ with $\sum_{i=1}^tk_i= k+1$,  the natural map $H^0(S,L) \to H^0(Z, L\otimes O_Z)$ is onto, where 
$(Z,O_Z)$ is a $0$-dimensional subscheme such that at  each point $z_i$:  $I_ZO_{S,z_i}$ is generated 
by $(x_i, y_i^{k_i})$, with $(x_i,y_i)$ local cordinates at $z_i$ on $S$.  Note that  $k =0, 1$ means respectively
 $L$  globally generated, $L$ very ample,(see \cite{BeltramettiFrancia}). In the sequel, we will need the following:
\begin{lem}
Let $S$ be a smooth complex projective surface and   $K_S$ be a canonical divisor on it.  Let  $H$ be a very ample divisor 
on $S$, such that $H^2 \geq 5$ and  $(S, O_S(H))$ does not contain lines and conics. Then we have the following  properties:
\begin{enumerate}
\item[(a)] the divisor $2H + K_S$ is very ample too;
\item[(b)] let $R \subset \vert 2H + K_S \vert $ be the locus of reducible curves, then $R$ is a closed subset of codimension $\geq 2$; 
\item[(c)] a general pencil $P \subset \vert 2H + K_S \vert $ has all irreducible elements and  the singular curves of $P$ have a unique node as singularities.
\end{enumerate}
\end{lem}
\begin{proof}
 Note that  since $H$ is very ample, to prove a),  it is enough that  $O_S(H + K_S)$ is a line bundle globally generated 
on $S$.  This is true for any pair $(S, O_S(H))$ which is not a  scroll or $({\bf P}^2, O_{{\bf P}^2}(i))$, $i=1,2$, see 
\cite{Sommese}. Let's examine  b). Let  $S^* \subset \vert 2H +K_S \vert$ be the locus of  singular curves, then $S^*$ 
is an irreducible variety and $codim S^* \geq 1$, see \cite{Harris2}, moreover  $R \subset S^*$, see \cite{Hartshorne} 
cor III 7.9. So property b) means either  $codim S^* \geq 2$ or   $R $ is a proper closed subset of $S^*$.
 Let $p \in S$ be any point,  let $\epsilon  \colon X \to S$ be the blow up of $S$ at the point $p$ with exceptional divisor
 $E$: assume  that  there  exists  a smooth irreducible curve in the linear system 
$\vert {\epsilon}^* (2H+K_S) -2 E \vert $, this  would  give us   an   irreducible curve having a unique node in $p$ 
in the linear system $\vert 2H+K_S  \vert $, which  implies   $R \not=  S^*$. For this it is enough to  request that
 ${\epsilon}^* (2H+K_S) -2 E$  is  ample and globally generated,  which is of course true if it is   very ample. This last
 property is achieved for every point $p$, whenever  $2H +K_S$ defines a $3-$rd order embedding, (i.e. it  is $3$ spanned),
 see \cite{Beltrametti}, prop.3.5. In particular, if $H^2 \geq 5$, $2H + K_S$ is $3-$spanned unless there exist an effective
 divisor $F$ on $S$ such that either H.F=1 and $F^2 = 0,-1,-2$ or $H.F=2$ and $F^2 =0$, see \cite{BeltramettiFrancia}.
 Since we assumed that there are no curves embedded by $H$ as lines or conics, this concludes b). Actually, we have also
 proved that a general element of $S^*$ is an irreducible curve with a unique node, which implies c).
\end{proof}
 \hfill\par
{\bf Remark:} Note that on  any surface $S$ we can easily find very ample line bundles satisfying the assumptions of the 
lemma: for any  very ample $H$, it is enough to take  $nH$ with $n\geq 3$. 
\section{Odd ramification coverings of smooth curves. }
 Let $X$ be a smooth, irreducible, complex projective curve of genus $g$. Let  $f \in {\bf C}(X)$ be  a non-constant 
meromorphic function on $X$  of degree $d$, then it defines a  holomorphic map $f \colon X \to  {\bf P}^1$,
 which is a ramified covering with branch locus $B \subset {\bf P}^1$ and ramification divisor $R \subset X$. 
Let $M(f)$ be the monodromy group of $f$, see 2.1. We say that $f$ is an  {\it odd ramification covering}  if all 
ramification points of $f$ have odd index. Note that if $f$ is an  odd ramification covering, then it has even monodromy,
 in fact all the generators of the group $M(f)$ can be decomposed in cycles of odd length. 
\hfill\par 
 \subsection{Constructing map with even monodromy}
 We  recall the  method used in  \cite{ArtebaniPirola} to produce  odd ramification  coverings. 
 A line bundle $L$ on $X$ is said a {\it spin bundle } if $L^{2}= K_X$, where $K_X$ denotes the canonical line bundle on
 $X$. Fix $3$ distinct points $p_1$, $p_2$, $p_3$ on $X$ and  define the divisor 
\begin{equation}
D = n_1p_1 + n_2p_2 + n_3p_3 \quad  n_i  \in {\bf N} \quad n_1 > n_2 > n_3 \geq 0;
\label{divisor}
\end{equation}
 set $d = deg D = n_1 + n_2 + n_3$ and   denote by $[D]$ 
 the  support of $D$, we have $deg [D] = k$ with  $k = 2$ or $3$. Let $L$ be a spin bundle on $X$: we consider the line
 bundle $L(D)$. Note that if $s $ is a global section in $H^0(X, L(D))$, then $s^2$ can be identified with a meromorphic
 form $\omega$ on $X$ having poles at the points of $[D]$. If  $\omega$ were an exact form, then there would be a 
non constant meromorphic function $f \in H^0(X, O_X(2D - [D]))$ on $X$, such that $\omega = df$. It is easy to verify
 that $f \colon X \to {\bf P}^1$ would be  a ramified covering with  odd ramification index at every point.
 Let us define set-theoretically
\begin{equation} {\cal A}(X,D)  = \{ s \in H^0(X, L(D)) \colon \ \ \ s^2 \ \ is \  \ exact \ \ \},
\end{equation}
\begin{equation} {\cal F}(X,D) = \{ f  \in {\bf C}(X) \colon \ \ \ df = s^2,   \ \ s \in  {\cal A}(X,D) \}.
\end{equation}
Note that ${\cal A}(X,D)$ is actually the zero scheme of the following map:
\label{psi}
\begin{equation} 
\psi \colon H^0(X,L(D)) \to H^1(X- [D], {\bf C})
\end{equation}  
sending each global section $s$ into the De Rham cohomology class $ [s^2 ]$ of the form $\omega = s^2$. Actually we will 
 consider the projectivization of ${\cal A}(X,D)$ 
\begin{equation} {\cal H}(X,D)  = \{ (s) \in  {\bf P}(H^0(X, L(D))) \colon \ \ \ s^2 \ \ is \  \ exact \ \ \}.
\end{equation}
We have the following results:
\begin{prop} Let $X$ be a smooth complex projective  curve of genus $g$, let $D$ be a divisor as in \ref{divisor} with degree $d$ 
and support of degree $k$. We assume that: $d > 8g + 3k -4$ and moreover 
 if $k= 2$  then  $2 n_i > 3g + 2 $  for $i=1,2$;   if $k = 3$ then  $2 n_i > 3g + 3 $, $i= 1,2,3$.
 Then ${\cal H}(X,D) \subset  {\bf P}^{d -1}$ is a complex
 projective  variety with the following properties:
\begin{enumerate}
\item[(i.)]   ${\cal H}(X,D)$ is an irreducible variety  of dimension $d -2g -k $ and its  singular locus
 $Sing({\cal H}(X,D))$ has dimension $h  < 4(g-1) + k$; 
\item[(ii.)]  ${\cal H}(X,D)\subset {\bf P}^{d -1}$ is a complete intersection of $2g + k -1$
linearly indipendent quadrics;
\item[(iii.)]  ${\cal H}(X,D)$ is a normal rationally connected variety.
\end{enumerate}
\end{prop}
\begin{proof}
Note that $\psi$ factors through the natural  linear map 
\begin{equation}
\theta \colon  {Sym}^2 H^0(X, L(D)) \to  H^1(X-[D],{\bf C}),
\end{equation}
defined as $\theta (s \otimes t) = [s.t]$, the De Rham cohomology class of $s.t$.  
 This implies that  ${\cal H}(X,D)$ is the zero locus of  homogeneous polinomials of degree $2$.
 Actually, ${\cal H}(X,D)$ can be seen as the zero locus of a global section $\sigma$ of the  following  
 vector bundle of rank $2g + k -1$ on ${\bf P}(H^0(X, L(D)))={\bf P}^{d -1}$:
\begin{equation}
E = H^1(X-[D],{\bf C}) \otimes O_ {{\bf P}^{d -1}}(2),
 \end{equation}
see \cite{Pirola}, pr.2.1. Note that the ideal sheaf ${\cal  I}_{ {\cal H}(X,D)} $ is the image of the dual map
 $\sigma^* \colon E^* \to O_{{\bf P}^{d -1}}$, hence it is locally generated by $2g + k -1$ elements.
 By studying the tangent map we can obtain  that, under the above assumptions,  
actually ${\cal H}(X,D)$ is irreducible of dimension $d -2g -k $ and moreover $dim Sing({\cal H}(X,D)) =  h < 4(g-1) + k$, 
see \cite{Pirola}, pr.5.1 and cor.5.3. This also implies that ${\cal H}(X,D)$ is a complete intersection of $ 2g + k -1$ 
quadrics  and concludes the proofs of i) and ii).
 ${\cal H}(X,D)$  is a normal variety since from i) it is regular in codimension 1, (see \cite{Hartshorne}, p.186).
 Finally, since ${\cal H}(X,D) \subset {\bf P}^{d -1}$ is an irreducible complete intersection of $2g + k -1$ quadrics,    by  pr. 2.2.1, it  is rationally connected if    we have:
$$ h  \leq d - 4g -2k,$$
this immediately follows from i), since we assumed $d > 8g + 3k -4$. 
 \end{proof}
 \hfill\par
Let $(s) \in {\cal H}(X,D)$: it defines a unique linear serie $g^1_m(s)$ on $X$ as follows:
$$ g^1_m (s)  = {\{ \lambda f + \mu \  = 0 \} }_{(\lambda, \mu) \in {\bf P}^1},$$
where $ f \in {\bf C}(X)$ and $d f = s^2$.
We have the following result:
\begin{prop} Let $X$ be a smooth complex projective  curve of genus $g$, let $D$ be a divisor with degree $d$ and support
  of degree $k$ as in \ref{divisor} . Assume that: $d > 6g + 2k -3$,   if $k=2$  then  $2n_i > 3g + 3$ for $i = 1,2$, if $k = 3$, 
then  $2n_i > 3g +4$ for $i= 1,2,3$,   moreover, the triple $(2n_1-1,2n_2-1,2n_3-1)$ is given by relatively 
prime integers. Then  for general  $(s)  \in {\cal H}(X,D)$ the linear serie $g^1_m(s)$ is base points free and defines 
an indecomposable  
finite morphism $F \colon X \to {\bf P}^1$ of degree $m = 2d -k$ with monodromy $M(F) = A_m$.
  \end{prop}
For the proof see \cite{ArtebaniPirola}, pr.3 and  th.1.
\section{Main constructions.}
In this section we will introduce some basic constructions, we will need in proving our main theorem.
\subsection{Lefschetz pencil.} 
Let $S$ be a smooth complex projective surface and let  $K_S$ be a canonical divisor on $S$. Let $H$
 be a very ample divisor on $S$ such that  $H^2 \geq 5$ and $(S,O_S(H))$ does not contain lines or conics. 
Set $g = p_a(2H + K_S)$ and  $N = (2H + K_S)^2 \geq H^2 \geq 5$. By lemma 2.4.1 c), we can   choose a  general pencil
 $P= { \{ C_t \} }_{t \in {\bf P}^1}$  in the linear system $\vert 2H + K_S \vert$,   with the following properties:
\begin{enumerate}
\item[(i.)] every curve in $P$ is irreducible;
\item[(ii.)] the generic curve in $P$ is a  smooth, irreducible,  complex projective curve of genus $g$;
\item[(iii.)] there are at most finitely many   singular curves  in $P$ and they have a unique node as singularities;
\item[(iv.)] every pair of curves $C_t$ and $C_{t'}$ of $P$ intersect trasversally, so that $P$ has $N$ distinct base 
points, $p_1$,...,$p_N$.
\end{enumerate}
We will call $P$ {\it  a Lefschetz pencil of irreducible curves  on $S$ of genus $g$}. Starting from these data $(S,H,P)$ 
we will introduce the following constructions. \hfill \par
\subsection{Construction 1} Let  $\hat S$ be the  smooth complex projective surface obtained by blowing up the base 
points of the pencil $P$:
\begin{equation}
\hat S= B_{p_1,p_2,...,p_N}(S).
\end{equation}
Let us   denote by $\epsilon \colon \hat S \to S$ the blow  up map, by $E_1$, ..,$E_N$ the exceptional curves, 
such that $E_i^2 = -1$ and $E_i \cdot E_j = 0$, for $i \not= j$, then $K_{\hat S} ={\epsilon}^* K_S  + E_1 +.. + E_N$.  
Note that the strict transforms of the curves of the pencil $P$ satisfies: $ \tilde{C_t}\cdot\tilde{C_{t'}}= 0,$ 
for any $t \not= t'$. Hence the pencil $P$ induces a surjective morphism 
\begin{equation}
\phi \colon \hat S \to {\bf P}^1,
\end{equation}
 with fibre $F_t = \tilde{C_t} $, for any $t \in  {\bf P}^1$,  $C_t \in P$. Moreover,  $\phi$ is actually a 
flat morphism and the exceptional curves $E_1,.. ,E_N$ in $\hat S$   turn out to be sections of the morphism $\phi$.
 We will define on $\hat S$ the line bundle  
\begin{equation}
L  = {\epsilon}^*O_S(H + K_S).
\end{equation}
Note that if  $F_t$ is any singular  fibre of $\phi$,   then its dualizing sheaf ${\omega}_{F_t}$ is a line bundle, 
since  we have ${\omega}_{F_t} = {({\omega}_{\hat S} + F_t)}_{\vert F_t}$, as for smooth fibres. 
It's easy to verify that for any fibre $F_t$ we have
\begin{equation}
L^2_{\vert F_t} \simeq {\omega}_{F_t},
\end{equation}
so we say  that  $L$ is  {\it a spin bundle relatively to  $\phi$}. 
We denote by $U \subset {\bf P}^1 $ the open subset corresponding to smooth fibres of $\phi$, 
set $\hat S_U = {\phi}^{-1}(U)$, then $\phi \colon \hat S_U \to U$ is a smooth morphism. 
We have proved the following \hfill \par
{\bf Claim 1}: {\it The smooth complex projective surface $\hat S$ is endowed with a surjective morphism
 $\phi \colon \hat S \to {\bf P}^1$, with smooth fibre $F_t$ of genus $g$,  
and a line bundle $L$ which is a spin bundle relatively to $\phi$.}
\subsection{Construction 2} Now let us  choose  three distinct exceptional curves $E_1$, $E_2$, $E_3$
 on the surface $\hat S$, and fix   integers $n_1 > n_2 > n_3 \geq 0$: we will   consider on $\hat S$ the line bundle
\begin{equation}
L(n_1E_1 + n_2E_2 + n_3E_3).
\end{equation}
Since each $E_i$ is  a  section of the morphism $\phi \colon \hat S \to {\bf P}^1$, for any  fibre $F_t$ of $\phi$ we have 
\begin{equation}
{L(n_1E_1 + n_2E_2 + n_3E_3)}_{\vert F_t} = L_t (D_t),
\end{equation}
where $D_t = n_1 p_1^t + n_2 p_2^t + n_3p_3^t$, with  $ p_i^t = {E_i}_{\vert F_t}$ and $ L_t  = L_{\vert F_t}$ is a 
spin bundle on $F_t$.  Set $d = deg (D_t) = n_1 + n_2 + n_3$ and  $[D_t] = $  the support of $D_t$, with 
$deg [D_t] = k$, $k = 2 $ or $3$.  We assume that: $d > 8g + 3k -4$,   if $k = 2$ then $2 n_i > 3g + 3$ for $i=1,2$; if $k = 3$ 
then  $2n_i > 3g + 4$ for $i= 1,2,3$, finally $(2n_1-1, 2n_2-1, 2n_3-1)$ are  relatively prime integers. 
For such $(d,k)$, for any smooth fibre $F_t$, by pr. 3.1.1,   we can  introduce the irreducible projective variety: 
\begin{equation}
{\cal H}(F_t,D_t) \subset  {\bf P}(H^0(F_t,L_t(D_t))) = {\bf P}^{d -1}_t.
\end{equation} 
{\bf Claim 2:} {\it There exists a complex  projective   variety ${\cal H}$ and a surjective   morphism 
$p \colon  {\cal H} \to {\bf P}^1$, with the following property:  let $ U \subset {\bf P}^1$ be the open subset 
corresponding to smooth fibres $F_t$ of $\phi$, for any $t \in U$, the  fibre  $p^{-1}(t)$ is  the projective variety
 ${\cal H}(F_t, D_t)$}.  \hfill \par
Let us  consider on $\hat S$ the line bundle  $L(n_1E_1 + n_2E_2 + n_3E_3)$ and  look at its  restriction
 $L_t (D_t)$ to any fibre $F_t$. Since $\forall t$ $F_t$ is irreducible and lies on a smooth surface, 
then $ deg (L_t(D_t))> 2p_a -2$, implies $h^1(F_t , L_t(D_t))=0$, see \cite{CataneseFranciosi}, so we can apply 
Riemann Roch theorem and obtain $h^0(F_t, L_t(D_t))= d$.  Since $\phi \colon \hat S \to {\bf P}^1$  is a flat morphism,
by Grauert 's theorem, (see \cite{Hartshorne}, p.288),  the sheaf 
\begin{equation}
{\cal F} = {\phi}_* ( L(n_1E_1 + n_2E_2 + n_3E_3))
 \end{equation}
is a locally free sheaf of rank $d$ on ${\bf P}^1$.  So we can introduce the associated projective space bundle
${\bf P}({\cal F})$ and the following smooth morphism
\begin{equation}
\pi \colon {\bf P}({\cal F})\to {\bf P}^1,
\end{equation}
whose fibre at $t$ is the projective space ${\bf P}(H^0(F_t,L_t(D_t))) = {\bf P}^{d -1}_t$. 
Let $O_{\cal F}(1)$ be the tautological line bundle on ${\bf P}({\cal F})$, i.e. 
$O_{\cal F}(1)_{\vert {\bf P}({\cal F}_t)}= O_{{\bf P}^{d -1}_t}(1)$. Let  $U \subset {\bf P}^1$ be the open subset 
where  $\phi$ is smooth and  $ \hat S_U = {\phi}^{-1}(U)$. Set   $W = \hat S_U- \{E_1, E_2, E_3\}$, 
 we can consider the restriction 
\begin{equation} 
\bar{\phi}= {\phi}_{\vert W} \colon W \to U,
\end{equation}
with fibre ${\bar{\phi}}^{-1}(t) = F_t - [D_t]$, for any $t \in U$. Since for any smooth fibre $F_t$, 
we have $h^1(F_t -[D_t],{\bf C}) = 2g + k -1$, the sheaf $R^1{\bar{\phi}}_* ({\bf C} )$ is actually  a vector bundle  on 
$U$ with  fibre $H^1(F_t -[D_t],{\bf C})$, set  
\begin{equation}
{\cal G}_U = R^1{\bar{\phi}}_* ({\bf C} ).
\end{equation}
Let $Sym^2 {{\cal F}}$ be the $2$ symmetric power of ${\cal F}$, we have the following natural maps:
\begin{equation}
{ \alpha} \colon {O_{\cal F}(-2)}_{\vert U} \to {Sym^2 {{\cal F}}}_{\vert U}  
\end{equation}
\begin{equation}
\Theta \colon {Sym^2 {{\cal F}}}_{\vert U}  \to {\cal G}_U,
\end{equation}
see  the  proof  of 3.1.1.
By composition we obtain   a non zero global section $\tau $ of the vector bundle  
${\cal G}_U \otimes O_{{\cal F}}(2)_{\vert U}$. We define the projective variety
\begin{equation}
{\cal H}_U  \subset   {{\bf P}({\cal F})}_{\vert U},
\end{equation}
as the zero locus of the  section $\tau$. 
It admits a  natural surjective  morphism $p_U  \colon {\cal H}_U \to U$,  whose fibre at $t$ is actually 
the projective variety ${\cal H}(F_t, D_t)$.  It' s easy to verify that  $p_U$ turns out to be  a proper flat morphism. 
Finally, let ${\cal H}$ be the scheme-theoretic closure of ${\cal H}_U$ into the projective variety  
${\bf P}({\cal F})$, then  ${\cal H}$ is a complex projective variety, moreover, since $U ={\bf P}^1 - \{ t_1,.., t_Q \}$,
 then   there exists  a flat morphism  $p \colon  {\cal H} \to {\bf P}^1$, which extends $p_U$,
 (see \cite{Hartshorne}, p.258).
\subsection{Claim 3: {\it  ${\cal H}$ admits a section $ \sigma$.}}
 Look at the surjective flat morphism $p \colon  {\cal H} \to {\bf P}^1$: for any  $t \in U$,
 the fibre $ p^{-1}(t) ={\cal H}(F_t, D_t)$  is a normal rationally connected variety, see pr. 3.1.1. This allows us 
to apply  theorem 2.2.2  to $p$ and to conclude that $p$ has a section, let us denote  it by $\sigma$,
\begin{equation} \sigma \colon  {\bf P}^1 \to {\cal H},
\end{equation}
with the following property:
for general $t \in U$,  the linear serie $g^1_m(t)$, defined by $\sigma(t)$, on the smooth  fibre $F_t$,  is base points 
free  of degree  $2d - k$. So, by   pr. 3.1.2, the associated map $F_t \to {\bf P^1}$ is indecomposable with   
 monodromy group $ A_{2d-k}$.
Note that under  the assumptions made in 4.3, $m$ is even and $m > 16g + 2$ if $k= 2$, while $m$ is odd and $m > 16g + 7$  
if $k = 3$.
\subsection{\bf Construction 3}{\it  There exist a smooth, complex, rational ruled surface $Y$  and  a  finite rational map  of degree $m = 2d -k$,  $\delta  \colon \hat S \to Y$ with the following property:   for a general smooth fibre $F_t$ the restriction ${\delta}_{\vert F_t}$ is given by  the linear serie $g^1_m(t)$ on $F_t$ and  
the following diagramm commutes}: 
\[
\begin{array}{c}
{\hat S} \\ {\downarrow} \\  {\bf  P}^1 \end{array}
\begin{array}{c}
\stackrel{\delta}{\rightarrow} \\   \\  \stackrel{id}{\rightarrow} \end{array}
\begin{array}{c}
{Y} \\ {\downarrow} \\  {\bf  P}^1 \end{array}
\]
{\it where the vertical arrows are respectively the morphism $\phi$ and the ruling $\pi$ of $Y$}.
\hfill \par
Since $\phi \colon \hat S_U \to U$ is a smooth morphism, we can consider the quasi projective variety $\hat S_U^{(m)}$ parametrizing the symmetric products $F_t^{(m)}$ of the smooth fibres $F_t$ of $\phi$.  There is natural map induced by $\phi$, which is  a smooth morphism 
\begin{equation}
{\phi}_m \colon \hat S_U^{(m)} \to U, 
\end{equation}
with smooth fibre $F_t^{(m)}$.   The  existence of  $\sigma$,  allows us to define a quasi projective variety  ${\cal I}$ as follows: 
\begin{equation}
{\cal I} = \{ (A, t ) \in \hat S_U^{(m)} \times U \colon  A \in g^1_m(t)  \}.
 \end{equation}
Let ${\pi}_1 \colon {\cal I} \to U$ the natural projection, then ${{\pi}_1}^{-1}(t) \simeq {\bf P}^1$ is the linear serie $g^1_m(t)$. So ${\cal I}$ is a quasi-projective surface endowed with a rational ruling ${\pi}_1$. Then let $Y$ be a smooth rational  ruled surface whose  ruling
\begin{equation}
{\pi} \colon Y \to {\bf P}^1,
\end{equation} 
restricts to $U$ is ${\pi}_1$, let  $F_t^Y $ denote the  fibre of ${\pi}$ at  $t$. Finally, we  define the rational map $\delta$: let $x \in \hat S_U$, then there exists a unique smooth fibre $F_t $ through $x$,  assume that $x$ is not a base point of the serie  $g^1_m(t)$, then  $\delta (x)$ is the unique divisor in $g^1_m(t)$ passing through the point $x$. It is easy to see that $\delta$ is a rational map.  Let $t \in  U $ be a general point, then  the fibre $F_t$ is smooth and the linear serie $g^1_m(t)$ is base points free of degree $m = 2d-k$, see 4.4.The restriction  ${\delta}_{\vert F_t} $  is actually the morphism associated to  $g^1_m(t)$:
\begin{equation}
{\delta}_{\vert F_t} \colon F_t \to  F_t^Y \simeq {\bf P}^1,
\end{equation}
so the map induced on the ${\bf P}^1$'s  must be the identity. Moreover, by  4.4,  for general $t \in U$,  the   monodromy group $M( {\delta}_{\vert F_t})$  is the alternating  group $ A_m$. 
\subsection{\bf Construction 4} The rational map $\delta \colon \hat S \to Y$  can be resolved with a finite number of blow ups as follows. Let $V \subset \hat S$ be an open subset where  $\delta$ is defined. Let ${\Gamma}_{\delta} \subset  \hat S \times Y $ be the closure of the graph of the morphism ${\delta}_{\vert V}$. ${\Gamma}_{\delta}$ is a projective variety, and it has two natural projections  ${\pi}_1 \colon  {\Gamma}_{\delta} \to \hat S$, which is a birational morphism, and ${\pi}_2 \colon  {\Gamma}_{\delta} \to Y$, which is a generically finite surjective morphism of degree $m$. Then there exist a smooth surface  $X $ and a  birational morphism $r \colon X \to {\Gamma}_{\delta}$ which is a resolution of singularities of  ${\Gamma}_{\delta}$, see \cite{Hironaka}. Hence  we have: 
\begin{enumerate}
\item[(i.)] $X$ is a smooth complex projective surface in the birational class of $S$;
\item[(ii.)] there exists a surjective morphism $\eta = \phi \cdot {\pi}_1 \cdot r \colon X \to {\bf P}^1$, whose  general smooth fibre is isomorphic to a general smooth fibre $F_t$ of $\phi$;
\item[(iii.)] there exists  a generically finite surjective  morphism, $f = {\pi}_2 \cdot r \colon X \to Y$, of degree $m$, such that  the restriction $f_{\vert F_t}$ is actually ${\delta}_{\vert F_t}$, for a general smooth fibre $F_t$. 
\end{enumerate}
So we have proved the following \hfill\par
{\bf Claim 4:} {\it We have a commutative diagramm: }
\[
\begin{array}{c}
{X} \\ {\downarrow} \\  {\bf  P}^1 \end{array}
\begin{array}{c}
 \stackrel{f}{\rightarrow} \\   \\  \stackrel{id}{\rightarrow} \end{array}
\begin{array}{c}
{Y} \\ {\downarrow} \\  {\bf  P}^1 \end{array}
\]
{\it where the vertical arrows are respectively the morphism $\eta$ and the ruling $\pi$ of $Y$,  such that for a general smooth fibre $F_t$, the  monodromy group $M( f_{\vert F_t})$ is the alternating  group $A_m$}. \hfill \par
\section{The main result.}
 \subsection{Technical lemma}
We start with a basic  lemma, which is  an  easy application of a topological result of Nori, (see \cite{Nori}, lemma 1.5). 
\begin{lem} Let $X$ be a smooth complex projective surface endowed with a surjective morphism  $\eta \colon X \to {\bf P}^1$  with general smooth fibre $F_t$. Let $Y$ be a smooth complex rational ruled  surface  with ruling ${\pi} \colon Y \to {\bf P}^1$, and fibre $F_t^Y \simeq {\bf P}^1$.  Assume that $f \colon X \to Y $ is a generically finite dominant morphism of degree $m$, such that the following diagramm commutes:
\[
\begin{array}{c}
{X} \\ {\downarrow} \\  {\bf  P}^1 \end{array}
\begin{array}{c}
 \stackrel{f}{\rightarrow} \\   \\  \stackrel{id}{\rightarrow} \end{array}
\begin{array}{c}
{Y} \\ {\downarrow} \\  {\bf  P}^1 \end{array}
\]
then the restriction of $f$ to general smooth fibres of $\eta$ completely induces the monodromy group $M(f)$, i.e.
$$ M(f,x) \simeq  M(f_{\vert F_t},x),$$
for a general smooth fibre $F_t$ and $x \in f(F_t)$, not a branch point of $f$.
\end{lem}
\begin{proof}
Let us consider the  morphism $f \colon X \to Y$, let  $R \subset X$ be the ramification divisor and  $B \subset Y$ the branch locus of $f$. The following map
\begin{equation}
q =  {f}_{\vert X- f^{-1}(B)}\colon X-f^{-1}(B) \to Y - B,
\end{equation}
is a covering of degree $m$ in the classic topology. Look at the restriction  to a general smooth fibre $F_t$ of $\eta$, by the above commutative diagramm, we have: 
$$ q_{\vert F_t} =  {f}_{\vert F_t - (F_t \cap f^{-1}(B))}\colon F_t - (F_t \cap f^{-1}(B)) \to F^Y_t - (B \cap F^Y_t),$$
since $R \cap F_t$ is actually the ramification divisor of ${f}_{\vert F_t}$ and $B \cap F^Y_t$ is the branch losus of ${f}_{\vert F_t}$, then $q_{\vert F_t}$ is a covering too of degree $m$. Now let's also restrict ${\pi}$  to $Y -B$:
\begin{equation}
{{\pi}}_{\vert Y - B} \colon   Y -B \to {\bf P}^1,
\end{equation}
by the above commutative diagramm, since the induced map on the ${\bf P}^1$'s  is the identity, $B$ cannot contain  a complete fibre. This allows us to conclude that the restriction ${{\pi}}_{\vert Y- B}$ is surjective too. Moreover, note that since ${\pi}$ is a ruling of a rational ruled surface, it admits a section: so it cannot have multiple fibres, that is every fibre must have a reduced component. So by   lemma (1.5), c) of \cite{Nori}, we have the following exact sequence between the fundamental groups:
\begin{equation}
{\pi}_1(F^Y_t - (B \cap F^Y_t)) \to {\pi}_1(Y - B) \to {\pi}_1({\bf P}^1),
\end{equation}
for a general smooth fibre $F^Y_t$. Since ${\pi}_1({\bf P}^1) = 0 $, this gives us a surjective map $s_t$: 
\begin{equation}
s_t \colon {\pi}_1(F^Y_t - (B \cap F^Y_t)) \to {\pi}_1(Y - B ).
\end{equation}
Let  $x \in Y - B$ be a  point such that $x \in  F^Y_t= f(F_t)$, for  a general smooth fibre $F_t$. We recall that the monodromy representation is the group homomorphism  
\begin{equation}
\rho(f,x) \colon {\pi}_1(Y - B ) \to  Aut (f^{-1}(x)),
\end{equation}
whose image is $M(f,x)$. The surjectivity of $s_t$ immediately implies
\begin{equation}
M(f,x) = M(f_{\vert F_t},x),
\end{equation}
for a general smooth fibre $F_t$ and for any $x \in f(F_t)$, $x \not\in B$. This concludes the proof.
\end{proof}
\subsection{Proof of Theorem 2. }
Let $m> 16g +7$ be any integer, we can find a pair of integers $(d,k)$ satisfying  the following properties:
$$
k = 2  \ or  \ 3,  \ \ \  d > 8g +3k -4, \   \ \  2d -k = m.$$
Since  $H^ 2 \geq 5$ and the pair  $(S, O_S(H))$ does not contain lines and conics, see 4.1, we can choose a Lefschetz pencil $P$ of irreducible curves of genus $g$, in the linear system $\vert 2H +K_S \vert $. We  can apply all constructions of sec. 4 to the data $(S, H, P)$, where $(d, k)$ are given as above. So we produce the following situation: $X$ is a  smooth complex projective surface,  birationally equivalent to $S$, endowed with a surjective morphism $\eta \colon X \to {\bf  P}^1 $, with smooth fibre $F_t$ of genus $g$, $Y$ is a smooth complex rational ruled surface,  with ruling ${{\pi}} \colon Y \to {\bf P}^1$, $f \colon X \to Y$ is a generically finite morphism  of degree $m = 2d-k$, finally  the following  diagramm commutes:  
\[
\begin{array}{c}
{X} \\ {\downarrow} \\  {\bf  P}^1 \end{array}
\begin{array}{c}
\stackrel{f}{\rightarrow} \\   \\  \stackrel{id}{\rightarrow} \end{array}
\begin{array}{c}
{Y} \\ {\downarrow} \\  {\bf  P}^1 \end{array}
\]
where the vertical maps are respectively $\eta$ and ${\pi}$. Moreover, for a general smooth fibre $F_t$ of $\eta$, the monodromy group $M(f_{\vert F_t})$ is the alternating group$ A_m$.
Note that all the assumptions of lemma 5.1.1 are verified, hence we have:
\begin{equation}
M(f,x) = M(f_{\vert F_t},x),
\end{equation} 
for a general smooth fibre $F_t$ and a point $x \in f(F_t)$, which is not a branch point.
Since  $M( f_{\vert F_t})= A_m$, for a general smooth fibre $F_t$, we can finally conclude that the monodromy group $M(f)$ is actually the alternating group $A_m$.
\hfill \par 
{\bf Remark:} Note that the above theorem works under the following more general hypothesis: let $H$ be an ample divisor, such that $2H+K_S$ is very ample and $2H+K_S$ defines a $3$-th order embedding, see lemma 2.4.1.\hfill \par  
\subsection{Surfaces of general type.}
  We would like to apply  the above result to surfaces of general type. Let $S$ be a minimal, smooth complex projective surface of general type with ample canonical divisor $K_S$. As it is well known, for some $n >0$  the pluricanonical map ${\phi}_{nK_S }$ is  an embedding; in order to apply theorem 2,   we will be interested in the smallest $n$ such  that  ${\phi}_{nK_S }$ is actually a $3$-th order embedding. In fact, in this   situation, if $n = 2t + 1 \geq 3$, we can find  a Lefschetz pencil $P$, of irreducible curves in the linear system $\vert n K_S \vert$, see 4.1, and apply our constructions of section 4 to the data $(S, O_S(tK_S), P)$.  At this hand,  we will use  the following result: 
\begin{lem} Let $S$ be a minimal surface of general type with ample cannonical divisor $K_S$.  
\begin{enumerate}
\item[i)] If $n \geq 5$, the divisor $n K_S$ is very ample, if  $p_g \geq 3$ and $K_S^2 \geq3$ then $3 K_S$ is very ample too;
\item[ii)] if $n \geq 5$ and $K_S^2 \geq3$, then $nK_S$ defines a $3$-th order embedding, moreover if $K_S^2 > 5$ then $3 K_S$ defines a $3$-th order embedding unless  there exists an effective divisor $F$ on $S$ such that  $K_S .F =2$ with $F^2 =0$.
\end{enumerate}
\end{lem} 
For the proofs see \cite{Catanese} for i) and \cite{BeltramettiFrancia} for ii).  \hfill\par
As an immediate consequence of our  result, we have the following:
\begin{teo}
Let $S$ be a minimal, smooth,  complex, connected,   projective surface of general type with ample canonical divisor 
$K_S$, with  $K_S^2 > 3$. Then for any $m > 16 (1 + 15 K_S^2) + 7$, there exist a smooth complex projective surface $X$, in the  birational class of $S$, and a  generically finite surjective  morphism, of degree $m$: $$f \colon X \to Y, $$
into a smooth complex rational ruled surface $Y$ such that the   monodromy group  $M(f)$ is the alternating group $A_m$. \hfill \par
Moreover, if $p_g \geq 3$ and $K_S^2 > 5$, and $S$ does not contain any effective divisor $F$ described in lemma 5.3.1,
 then $m > 16 (1 + 6 K_S^2) + 7$.    
\end{teo}

\noindent {\sc Sonia Brivio and Gian Pietro Pirola }\ \ Dipartimento di Matematica, Universit\`a di Pavia\ \ 
via Ferrata 1, 27100 Pavia, Italia \ \ e-mail: brivio@dimat.unipv.it; pirola@dimat.unipv.it 

\end{document}